\documentclass[12pt,reqno]{amsart}
\usepackage{amsmath, amsthm, amscd, amsfonts, amssymb, graphicx, color}
\usepackage[bookmarksnumbered, colorlinks, plainpages]{hyperref}

\theoremstyle{definition}

\theoremstyle{remark}

\numberwithin{equation}{section}

\begin{document}
\begin{center}
{\textbf{$*$-Conformal $\eta$-Ricci Soliton on Sasakian manifold}}
\end{center}
\vskip 0.3cm
\begin{center}By\end{center}\vskip 0.3cm
\begin{center}
{Soumendu Roy \footnote{The first author is the corresponding author.}, Santu Dey, Arindam~~Bhattacharyya and Shyamal Kumar Hui}
\end{center}
\vskip 0.3cm
\address[Soumendu Roy]{Department of Mathematics,Jadavpur University, Kolkata-700032, India}
\email{soumendu1103mtma@gmail.com}

\address[Santu Dey]{Department of Mathematics,Techno Main Saltlake, EM 4/1,Sector- V, Saltlake,
Kolkata-700091, India}
\email{santu.mathju@gmail.com}

\address[Arindam Bhattacharyya]{Department of Mathematics,Jadavpur University, Kolkata-700032, India}
\email{bhattachar1968@yahoo.co.in}

\address[Shyamal Kumar Hui] {Department of Mathematics, The University of Burdwan, Burdwan, 713104, India.}
\email{skhui@math.buruniv.ac.in}
\vskip 0.5cm
\begin{center}
\textbf{Abstract}\end{center}
In this paper we study $*$-Conformal $\eta$-Ricci soliton on Sasakian manifolds. Here, we discuss some curvature properties on Sasakian manifold admitting $*$-Conformal $\eta$-Ricci soliton. We obtain some significant results on $*$-Conformal $\eta$-Ricci soliton in Sasakian manifolds satisfying $R(\xi, X)\cdot S = 0$, $S(\xi, X)\cdot R = 0$, $\overline{P}(\xi, X)\cdot S$ $=0$, where  $\overline{P}$ is Pseudo-projective curvature tensor.The conditions for $*$-Conformal $\eta$-Ricci soliton on $\phi$-conharmonically flat and $\phi$-projectively flat Sasakian manifolds have been obtained in this article. Lastly we give an example of 5-dimensional Sasakian manifolds satisfying $*$-Conformal $\eta$-Ricci soliton.$*$-Conformal $\eta$-Ricci soliton, $\eta$-Einstein manifold, Sasakian manifold.  
\vspace {0.3cm}
\section{\textbf{Introduction}}
The notion of Ricci flow was first introduced by Hamilton \cite{rsham} in 1982.The Ricci flow is an evolution equation for metrics
on a Riemannian manifold. The Ricci flow equation is given by\\
\begin{equation}
\frac{\partial g}{\partial t} = -2S
\end{equation}\\
on a compact Riemannian manifold $M$ with Riemannian metric $g$.\\\\
A self-similar solution to the Ricci flow \cite{rsham}, \cite{topping} is called a Ricci soliton \cite{rsha} if it moves only by a one parameter family of diffeomorphism and scaling. The Ricci soliton equation is given by\\
\begin{equation}
\pounds_V g + 2S + 2\lambda g=0,
\end{equation}\\
where $\pounds_V$ is the Lie derivative in the direction of $V$, $S$ is Ricci tensor, $g$ is Riemannian metric, $V$ is a vector field and $\lambda$ is a scalar. The Ricci soliton is said to be shrinking, steady and expanding accordingly as $\lambda$ is negative, zero and positive respectively.\\\\
A.E. Fischer during 2003-2004 developed the concept of conformal Ricci flow \cite{aefi}  which is a variation of the
classical Ricci flow equation that modifies the unit volume constraint of that equation to a scalar curvature constraint. The conformal Ricci flow on $M$ is defined by the equation \cite{aefi}\\
\begin{equation}
\frac{\partial g}{\partial t} + 2(S + \frac{g}{n}) = -pg
\end{equation}\\
and $r(g) = -1,$\\
where $M$ is considered as a smooth closed connected oriented $n$-manifold, $p$ is a scalar non-dynamical field(time dependent scalar field), $r(g)$ is the scalar curvature of the manifold and $n$ is the dimension of manifold.\\\\
The notion of Conformal Ricci soliton equation was introduced by N. Basu and A. Bhattacharyya \cite{nbab} in 2015 and the equation is given by \\
\begin{equation}
\pounds_V g + 2S = [2\lambda - (p + \frac{2}{n})]g,
\end{equation}\\
 where $\lambda$ is constant. \\
The equation is the generalization of the Ricci soliton equation and it also satisfies the conformal Ricci flow equation.\\\\
In 2009, Jong Taek Cho and Makoto Kimura introduced the notion of $\eta$-Ricci soliton \cite{joma}, given by the equation:
\begin{equation}
  \pounds_\xi g + 2S =2\lambda g+2 \mu \eta \otimes \eta,
\end{equation}
for constants $\lambda$ and $\mu$.\\
In 2018, Mohd Danish Siddiqi \cite{mohd} introduced the notion of Conformal $\eta$- Ricci soliton as:
\begin{equation}
  \pounds_\xi g + 2S +[2\lambda - (p + \frac{2}{n})]g+2 \mu \eta \otimes \eta=0,
\end{equation}
where  $\pounds_\xi$ is the Lie derivative along the vector field $\xi$ , $S$ is the Ricci tensor, $\lambda$, $\mu$ are contants, $p$ is a scalar non-dynamical field(time dependent scalar field)and $n$ is the dimension of manifold.\\
The notion of $*$-Ricci tensor on almost Hermitian manifolds and $*$-Ricci tensor of real hypersurfaces in non-flat complex space were introduced by Tachibana \cite{tachi} and Hamada \cite{hama} respectively where the $*$-Ricci tensor is defined by:
\begin{equation}
  S^*(X,Y)=\frac{1}{2}(trace \{\phi\circ R(X,\phi Y)\}),
\end{equation}
for all vector fields $X,Y$ on $M^n$ and $\phi$ is a (1,1)-tensor field.\\
If $S^*(X,Y)=\lambda g(X,Y)+\mu \eta(X)\eta(Y)$ for all vector fields $X,Y$ and $\lambda$, $\mu$ are smooth functions, then the manifold is called $*-\eta$-Einstein manifold.\\
Further if $\mu=0$ i.e $S^*(X,Y)=\lambda g(X,Y)$ for all vector fields $X,Y$ then the manifold becomes $*$-Einstein.\\
In 2014 Kaimakamis and Panagiotidou \cite{kaipan} introduced the notion of $*$-Ricci soliton which can be defined as:
\begin{equation}
  \pounds_V g + 2S^* + 2\lambda g=0,
\end{equation}
for all vector fields $X,Y$ on $M^n$ and $\lambda$ being a constatnt.\\\\
Now we define the notion of $*$-Conformal $\eta$ -Ricci soliton as:
\begin{equation}
  \pounds_\xi g + 2S^* +[2\lambda - (p + \frac{2}{n})]g+2 \mu \eta \otimes \eta=0,
\end{equation}
where  $\pounds_\xi$ is the Lie derivative along the vector field $\xi$, $S^*$ is the $*$- Ricci tensor and $\lambda$, $\mu$, $p$, $n$ are as  defined in (1.6).\\
The Riemannian-Christoffel curvature tensor $R$ \cite{ozgu}, the conharmonic
curvature tensor $H$ \cite{ishii}, the projective curvature tensor $P$ \cite{yano} and the Pseudo-projective curvature tensor $\overline{P}$ \cite{prasad} in a Riemannian manifold $(M, g)$ are defined by:
\begin{equation}
  R(X, Y )Z = \nabla_X\nabla_Y Z - \nabla_Y \nabla_X Z - \nabla_{[X,Y]} Z,
\end{equation}
\begin{eqnarray}
 H(X, Y )Z &=& R(X, Y )Z -\frac{1}{(n-2)}[g(Y,Z)QX - g(X,Z)QY\nonumber \\
           &+& S(Y,Z)X - S(X,Z)Y],
\end{eqnarray}
 \begin{equation}
       P(X, Y )Z = R(X, Y )Z -\frac{1}{(n-1)}[g(Y,Z)QX - g(X,Z)QY ],
     \end{equation}
     and \begin{eqnarray}
            \overline{P}(X,Y))Z &=& aR(X,Y)Z+b[S(Y,Z)X-S(X,Z)Y] \nonumber \\
                                &-& \frac{r}{n}(\frac{a}{n-1}+b)[g(Y,Z)X-g(X,Z)Y],
         \end{eqnarray}
where $Q$ is the Ricci operator, defined by $S(X, Y ) = g(QX, Y )$, $S$ is
the Ricci tensor, $r = tr(S)$ is the scalar curvature, where $tr(S)$ is the trace of $S$, $a, b \neq 0$ are constants and $X, Y, Z \in \chi(M)$, $\chi(M)$ being
the Lie algebra of vector fields of M.\\\\
The paper is organized as follows: After introduction, section 2 consists of basic definitions of Sasakian manifolds. Section 3 devotes $*$-Conformal $\eta$-Ricci soliton on Sasakian manifolds. In section 4, we establish some significant results on $*$-Conformal $\eta$-Ricci soliton in Sasakian manifolds satisfying $R(\xi, X)\cdot S = 0$, $S(\xi, X)\cdot R = 0$, $\overline{P}(\xi, X)\cdot S = 0$, where  $\overline{P}$ is Pseudo-projective curvature tensor. Moreover it is shown that the Sasakian manifold, admitting $*$-Conformal $\eta$-Ricci soliton becomes $\eta$-Einstein manifold when it satisfies $S(\xi, X)\cdot R = 0$ and Einstein manifold when it satisfies $R(\xi, X)\cdot S = 0$ and $\overline{P}(\xi, X)\cdot S = 0$. We also obtain the conditions for $*$-Conformal $\eta$-Ricci soliton on $\phi$-conharmonically flat and $\phi$-projectively flat Sasakian manifolds. It is also shown that the Sasakian manifold, admitting $*$-Conformal $\eta$-Ricci soliton becomes $\eta$-Einstein manifold when it is $\phi$-conharmonically flat and Einstein manifold when it is $\phi$-projectively flat. In last section, we give an example of 5-dimensional Sasakian manifolds satisfying $*$-Conformal $\eta$-Ricci soliton.\\\\
\vspace {0.3cm}
\section{\textbf{Preliminaries}}
Let $M$ be a (2$n$+1)dimensional connected almost contact metric manifold with an almost contact metric structure $(\phi, \xi, \eta, g)$ where $\phi$ is a $(1,1)$ tensor field, $\xi$ is a vector field, $\eta$ is a 1-form  and $g$ is the compatible Riemannian metric such that\\
\begin{equation}
\phi^2(X) = -X + \eta(X)\xi, \eta(\xi) = 1, \eta \circ \phi = 0, \phi \xi = 0,
\end{equation}\\
\begin{equation}
g(\phi X,\phi Y) = g(X,Y) - \eta(X)\eta(Y),
\end{equation}\\
\begin{equation}
g(X,\phi Y) = -g(\phi X,Y),
\end{equation}\\
\begin{equation}
g(X,\xi) = \eta(X),
\end{equation}\\
for all vector fields $X, Y \in \chi(M).$\\\\
A (2$n$+1) dimensional almost contact manifold $M$ with $(\phi, \xi, \eta, g)$ structure is said to be a contact metric manifold iff $d\eta(X,Y)=g(X,\phi Y)$.\\
Also an almost contact metric structure of $M$ is said to be normal if
\begin{equation}
  2d\eta(X,Y)\xi+[\phi,\phi](X,Y)=0,
\end{equation}
where $[\phi,\phi]$ is Nijenhuis tensor.\\
A normal contact metric manifold is said to be a Sasakian manifold.\\
In a (2$n$+1) dimensional Sasakian manifold the following conditions hold:
\begin{equation}
  \nabla_X \xi=-\phi X,
\end{equation}
\begin{equation}
  R(X,Y)Z=g(Y,Z)X-g(X,Z)Y,
\end{equation}
\begin{equation}
  R(X,Y)\xi=\eta(Y)X-\eta(X)Y,
\end{equation}
\begin{equation}
  R(\xi,X)Y=g(X,Y)\xi-\eta(Y)X,
\end{equation}
\begin{equation}
  \eta(R(X,Y)Z)=g(Y,Z)\eta(X)-g(X,Z)\eta(Y),
\end{equation}
\begin{equation}
  (\nabla_X \eta)Y=-g(\phi X,Y),
\end{equation}
where $R$ is the Riemannian curvature tensor and $\nabla$ is the Levi-Civita  connection associated with $g$.\\
As every Sasakian manifold is K-contact then we have in a Sasakian manifold,
\begin{equation}
  (\pounds_\xi g)(X,Y)=0,
\end{equation}
where $\pounds_\xi$ is the Lie derivative along the vector field $\xi$ which is a Killing vector field.\\
\vspace {0.3cm}
\section{\textbf{$*$-Conformal $\eta$-Ricci soliton on Sasakian manifold}}
Let $ (M,\phi,\xi,\eta,g)$ be a (2$n$+1) dimensional Sasakian manifold.Consider the  $*$-Conformal $\eta$-Ricci soliton on $M$ as:
\begin{equation}
  \pounds_\xi g + 2S^* +[2\lambda - (p + \frac{2}{2n+1})]g+2 \mu \eta \otimes \eta=0,
\end{equation}
Then from (2.12), we get,
\begin{equation}
  S^*(X,Y)=-[\lambda-\frac{1}{2}(p + \frac{2}{2n+1})]g(X,Y)-\mu \eta(X)\eta(Y).
\end{equation}
In paper \cite{amdh}, Lemma 3.1, authors have proved that the $*$-Ricci tensor on a (2$n$+1)-dimensional Sasakian manifold  $ (M,\phi,\xi,\eta,g)$ is given by,
\begin{equation}
  S^*(X,Y)=S(X,Y)-(2n-1)g(X,Y)-\eta(X)\eta(Y),
\end{equation}
for all vector fields $X, Y$ on $M$.\\
Using this relation (3.3) we can write (3.2) as,
\begin{equation}
  S(X,Y)=[2n-1-\lambda+\frac{1}{2}(p + \frac{2}{2n+1})]g(X,Y)-(\mu-1)\eta(X)\eta(Y).
\end{equation}
From this equation (3.4), we get,
\begin{equation}
  S(X,\xi)=[2n-\lambda-\mu+\frac{1}{2}(p + \frac{2}{2n+1})]\eta(X).
\end{equation}
Again we also know in a (2$n$+1) dimensional Sasakian manifold the Ricci tensor field satisfies:
\begin{equation}
  S(X,\xi)=2n\eta(X).
\end{equation}
Then from (3.5) and (3.6), we get,
\begin{equation}
  \lambda+\mu=\frac{1}{2}(p + \frac{2}{2n+1}).
\end{equation}
Now we know,
\begin{equation}
  (\nabla_X S)(Y,Z)=XS(Y,Z)-S(\nabla_X Y,Z)-S(Y,\nabla_X Z),
\end{equation}
for any $X, Y, Z \in \chi(M)$.\\
Then replacing the expression of S from (3.4), we get,
\begin{equation}
   (\nabla_X S)(Y,Z)=-(\mu-1)[\eta(Y)(\nabla_X \eta)Z+\eta(Z)(\nabla_X \eta)Y].
\end{equation}
Using (2.11), we get,
\begin{equation}
   (\nabla_X S)(Y,Z)=(\mu-1)[\eta(Y)g(\phi X,Z)+\eta(Z)g(\phi X,Y)].
\end{equation}
Similarly we can obtain,
\begin{equation}
   (\nabla_Y S)(Z,X)=(\mu-1)[\eta(Z)g(\phi Y,X)+\eta(X)g(\phi Y,Z)],
\end{equation}
and
\begin{equation}
   (\nabla_Z S)(X,Y)=(\mu-1)[\eta(X)g(\phi Z,Y)+\eta(Y)g(\phi Z,X)].
\end{equation}
Then adding (3.10), (3.11), (3.12) and using (2.3) we get,
\begin{equation}
  (\nabla_X S)(Y,Z)+ (\nabla_Y S)(Z,X)+ (\nabla_Z S)(X,Y)=0,
\end{equation}
for any $X, Y, Z \in \chi(M)$.\\
Thus we can state the following theorem:\\\\
\textbf{Theorem 3.1.} {\em If a (2$n$+1) dimensional Sasakian manifold $ (M,\phi,\xi,\eta,g)$ admits $*$-Conformal $\eta$-Ricci Soliton
then the manifold has cyclic Ricci tensor i.e. $  (\nabla_X S)(Y,Z)+ (\nabla_Y S)(Z,X)+ (\nabla_Z S)(X,Y)=0$, for any $X,Y,Z \in \chi(M)$.}\\\\
Now if $\nabla S =0$, then taking $Z=\xi$ in the expression of $\nabla S$ from (3.10), we obtain,
\begin{equation}
  (\mu-1)g(\phi X,Y)=0.
\end{equation}
$\forall X, Y \in \chi(M)$.\\\\
Hence $\mu=1$.\\\\
Then from (3.7), we obtain,
 $\lambda=\frac{1}{2}(p + \frac{2}{2n+1})-1$.\\\\
 This leads to the following:\\\\
 \textbf{Theorem 3.2.} {\em If a (2$n$+1) dimensional Sasakian manifold $ (M,\phi,\xi,\eta,g)$ is Ricci symmetric i.e $\nabla S=0$ and admits $*$-Conformal $\eta$-Ricci Soliton then $\mu=1$ and  $\lambda=\frac{1}{2}(p + \frac{2}{2n+1})-1$.\\\\}
 \section{\textbf{Curvature properties on Sasakian manifold satisfying\\ $*$-Conformal $\eta$-Ricci soliton }}
The condition $R(\xi,X)\cdot S=0$ implies that,
\begin{equation}
  S(R(\xi,X)Y,Z)+S(Y,R(\xi,X)Z)=0,
\end{equation}
for any $X, Y, Z \in \chi(M)$.\\
Replacing the expression of $S$ from (3.4) and from the symmetries of $R$ we get,
\begin{equation}
  (\mu-1)[g(X,Y)\eta(Z)+g(X,Z)\eta(Y)-2\eta(X)\eta(Y)\eta(Z)]=0,
\end{equation}
for any $X, Y, Z \in \chi(M)$.\\
For $Z=\xi$ we have,
\begin{equation}
  (\mu-1)[g(X,Y)-\eta(X)\eta(Y)]=0,
\end{equation}
for any $X, Y\in \chi(M)$.\\
Then using (2.2) we have,
\begin{equation}
  (\mu-1)g(\phi X,\phi Y)=0,
\end{equation}
for any $X, Y\in \chi(M)$.\\
It follows that $(\mu-1)=0$ i.e,
\begin{equation}
 \mu=1.
\end{equation}
Then from (3.7) we get,
\begin{equation}
  \lambda=\frac{1}{2}(p + \frac{2}{2n+1})-1.
\end{equation}
So we can state the following theorem:\\\\
\textbf{Theorem 4.1.} {\em If a (2$n$+1) dimensional Sasakian manifold $ (M,\phi,\xi,\eta,g)$ satisfies $R(\xi,X)\cdot S=0$ and admits  $*$-Conformal $\eta$-Ricci Soliton then $\mu=1$ and $\lambda=\frac{1}{2}(p + \frac{2}{2n+1})-1$.\\\\}
Now from (4.5), (4.6) and (3.4) we obtain,
\begin{equation}
  S(X,Y)=2ng(X,Y).
\end{equation}
for any $X, Y \in \chi(M)$.\\
Then we have,\\
\textbf{Corollary 4.2.} {\em If a (2$n$+1) dimensional Sasakian manifold $ (M,\phi,\xi,\eta,g)$ satisfies $R(\xi,X)\cdot S=0$ and admits  $*$-Conformal $\eta$-Ricci Soliton then the manifold becomes Einstein manifold.\\\\}
Again the condition $S(\xi,X)\cdot R=0$ implies that,
\begin{multline}
  S(X,R(Y,Z)W)\xi-S(\xi,R(Y,Z)W)X+S(X,Y)R(\xi,Z)W-S(\xi,Y)R(X,Z)W \\
  +S(X,Z)R(Y,\xi)W-S(\xi,Z)R(Y,X)W+S(X,W)R(Y,Z)\xi-S(\xi,W)R(Y,Z)X\\
  =0,
\end{multline}
for any $X, Y, Z, W \in \chi(M)$.\\
Taking the inner product with $\xi$,the above equation becomes,
\begin{multline}
  S(X,R(Y,Z)W)-S(\xi,R(Y,Z)W)\eta(X)+S(X,Y)\eta(R(\xi,Z)W)\\
  -S(\xi,Y)\eta(R(X,Z)W)+S(X,Z)\eta(R(Y,\xi)W)-S(\xi,Z)\eta(R(Y,X)W)\\
  +S(X,W)\eta(R(Y,Z)\xi)-S(\xi,W)\eta(R(Y,Z)X)=0,
\end{multline}
for any $X, Y, Z, W \in \chi(M)$.\\
Replacing the expression of $S$ from (3.4) and taking $Z=\xi$, $W=\xi$ we get,
\begin{equation}
  [4n-1-2\lambda-\mu+(p + \frac{2}{2n+1})][g(X,Y)-\eta(X)\eta(Y)] = 0,
\end{equation}
for any $X, Y \in \chi(M)$.\\
Then using (2.2) we have,
\begin{equation}
  [4n-1-2\lambda-\mu+(p + \frac{2}{2n+1})]g(\phi X,\phi Y) = 0,
\end{equation}
for any $X, Y \in \chi(M)$.\\
Now using (3.7) the above equation becomes,
\begin{equation}
 [4n-1-\lambda+\frac{1}{2}(p + \frac{2}{2n+1})]g(\phi X,\phi Y)=0,
\end{equation}
for any $X, Y \in \chi(M)$.\\
It follows that $4n-1-\lambda+\frac{1}{2}(p + \frac{2}{2n+1})=0$ i.e,
\begin{equation}
  \lambda=4n+\frac{1}{2}(p + \frac{2}{2n+1})-1.
\end{equation}
and then from (3.7), we obtain,
\begin{equation}
  \mu=1-4n.
\end{equation}
This leads to the following:\\\\
\textbf{Theorem 4.3.} {\em If a (2$n$+1) dimensional Sasakian manifold $ (M,\phi,\xi,\eta,g)$ satisfies $S(\xi,X)\cdot R=0$ and admits  $*$-Conformal $\eta$-Ricci Soliton then $\mu=1-4n$ and $\lambda=4n+\frac{1}{2}(p + \frac{2}{2n+1})-1$.\\\\}
Now from (4.13), (4.14) and (3.4) we obtain,
\begin{equation}
  S(X,Y)=-2ng(X,Y)+4n\eta(X)\eta(Y),
\end{equation}
for any $X, Y \in \chi(M)$.\\
Then we have,\\\\
\textbf{Corollary 4.4.} {\em If a (2$n$+1) dimensional Sasakian manifold $ (M,\phi,\xi,\eta,g)$ satisfies $S(\xi,X)\cdot R=0$ and admits  $*$-Conformal $\eta$-Ricci Soliton then the manifold becomes $\eta$-Einstein manifold.\\\\}
The condition $\overline{P}(\xi,X)\cdot S=0$ implies that\\
\begin{equation}
  S(\overline{P}(\xi,X)Y,Z)+S(Y,\overline{P}(\xi,X)Z)=0.
\end{equation}
for any $X, Y, Z \in \chi(M)$, where $\overline{P}$ is the Pseudo-projective curvature tensor in $M$.\\
Now in a (2$n$+1) dimensional Sasakian manifold,
\begin{eqnarray}
   \overline{P}(\xi,X)Y &=& aR(\xi,X)Y+b[S(X,Y)\xi-S(\xi,Y)X] \nonumber \\
                        &-& \frac{r}{2n+1}(\frac{a}{2n}+b)[g(X,Y)\xi-g(\xi,Y)X],
\end{eqnarray}
where $a, b\neq 0$ are constants.\\
Then using (2.9), (3.4) and (3.5) we get,
\begin{eqnarray}
   \overline{P}(\xi,X)Y &=& [a-\frac{r}{2n+1}(\frac{a}{2n}+b)][g(X,Y)\xi-\eta(Y)X]\nonumber\\
                        &+& b[ 2n-1-\lambda+\frac{1}{2}(p+\frac{2}{2n+1})][g(X,Y)\xi-\eta(Y)X]\nonumber\\
                        &+& b(\mu-1)[\eta(Y)X-\eta(X)\eta(Y)\xi].
\end{eqnarray}
Replacing the expression of $S$ from (3.4) and using (4.18), we obtain,
\begin{multline}
  S(\overline{P}(\xi,X)Y,Z)  \\
  =[2n-1-\lambda+\frac{1}{2}(p+\frac{2}{2n+1})][a-\frac{r}{2n+1}(\frac{a}{2n}+b)][g(X,Y)\eta(Z)-g(X,Z)\eta(Y)]\\
  +b[2n-1-\lambda+\frac{1}{2}(p+\frac{2}{2n+1})]^2[g(X,Y)\eta(Z)-g(X,Z)\eta(Y)]\\
  -b(\mu-1)[2n-1-\lambda+\frac{1}{2}(p+\frac{2}{2n+1})][g(X,Y)\eta(Z)-g(X,Z)\eta(Y)]\\
  -(\mu-1)[a-\frac{r}{2n+1}(\frac{a}{2n}+b)][g(X,Y)\eta(Z)-\eta(X)\eta(Y)\eta(Z)].\\
\end{multline}
for any $X, Y, Z \in \chi(M)$.\\
As from (3.4) we have $S(X,Y)=S(Y,X)$ for any $X,Y\in\chi(M)$, then we can write,
\begin{equation}
  S(Y,\overline{P}(\xi,X)Z)=S(\overline{P}(\xi,X)Z,Y)
\end{equation}
for any $X, Y, Z \in \chi(M)$.\\
Then similarly from (4.19) by interchanging $Y,Z$, we get,
\begin{multline}
  S(\overline{P}(\xi,X)Z,Y)  \\
  =[2n-1-\lambda+\frac{1}{2}(p+\frac{2}{2n+1})][a-\frac{r}{2n+1}(\frac{a}{2n}+b)][g(X,Z)\eta(Y)-g(X,Y)\eta(Z)]\\
  +b[2n-1-\lambda+\frac{1}{2}(p+\frac{2}{2n+1})]^2[g(X,Z)\eta(Y)-g(X,Y)\eta(Z)]\\
  -b(\mu-1)[2n-1-\lambda+\frac{1}{2}(p+\frac{2}{2n+1})][g(X,Z)\eta(Y)-g(X,Y)\eta(Z)]\\
  -(\mu-1)[a-\frac{r}{2n+1}(\frac{a}{2n}+b)][g(X,Z)\eta(Y)-\eta(X)\eta(Y)\eta(Z)].\\
\end{multline}
for any $X, Y, Z \in \chi(M)$.\\
Now using (4.19), (4.20) and (4.21), (4.16) becomes,
\begin{multline}
  (\mu-1)[a-\frac{r}{2n+1}(\frac{a}{2n}+b)][g(X,Y)\eta(Z)+g(X,Z)\eta(Y)-2\eta(X)\eta(Y)\eta(Z)]=0.\\
\end{multline}
for any $X, Y, Z \in \chi(M)$.\\
Taking $Z=\xi$ in the above equation, we get,
\begin{equation}
  (\mu-1)[a-\frac{r}{2n+1}(\frac{a}{2n}+b)][g(X,Y)-\eta(X)\eta(Y)]=0.
\end{equation}
for any $X, Y \in \chi(M)$.\\
Then by using (2.2) we obtain,
\begin{equation}
  (\mu-1)[a-\frac{r}{2n+1}(\frac{a}{2n}+b)]g(\phi X,\phi Y)=0.
\end{equation}
for any $X, Y \in \chi(M)$.\\\\
Hence we get, $(\mu-1)[a-\frac{r}{2n+1}(\frac{a}{2n}+b)]=0.$\\\\
Then either $\mu-1=0$,or $a-\frac{r}{2n+1}(\frac{a}{2n}+b)=0.$\\\\
Which implies that, either,
\begin{equation}
  \mu=1.
\end{equation}
or,
\begin{equation}
  r=\frac{2n(2n+1)a}{a+2nb}.
\end{equation}
As in (1.3), $r = -1$, then from (4.26), we get,
\begin{equation}
  [2n(2n+1)+1]a + 2nb = 0.
\end{equation}
Now if $\mu=1$ then from (3.7) we get,
\begin{equation}
\lambda=\frac{1}{2}(p + \frac{2}{2n+1})-1.
\end{equation}
Then we can state the following theorem:\\\\
\textbf{Theorem 4.5.} {\em If a (2$n$+1) dimensional Sasakian manifold $ (M,\phi,\xi,\eta,g)$ satisfies $\overline{P}(\xi,X)\cdot S=0$ and admits  $*$-Conformal $\eta$-Ricci Soliton then either $\mu=1$ and $\lambda=\frac{1}{2}(p + \frac{2}{2n+1})-1$ or $[2n(2n+1)+1]a + 2nb = 0$, where $a,b\neq 0$ are constants.}\\\\
Now if $\mu=1$ and $\lambda=\frac{1}{2}(p + \frac{2}{2n+1})-1$ then from (3.4) we get,
\begin{equation}
S(X,Y)=2ng(X,Y).
\end{equation}
So we have,\\
\textbf{Note:} {\em If a (2$n$+1) dimensional Sasakian manifold $ (M,\phi,\xi,\eta,g)$ satisfies $\overline{P}(\xi,X)\cdot S=0$ and admits  $*$-Conformal $\eta$-Ricci Soliton then the manifold becomes Einstein manifold, provided $[2n(2n+1)+1]a + 2nb \neq 0$, where $a, b\neq 0$ are constants.}\\\\
\textbf{Definition 4.6.} A differentiable manifold $(M^n,g), n>3$, satisfying the condition
\begin{eqnarray}
  \phi^2H(\phi X,\phi Y)\phi Z &=& 0 \nonumber
\end{eqnarray}
is called $\phi$-conharmonically flat \cite{singh}, where $H$ is the conharmonic curvature tensor of $M^n$.\\
Let $ (M,\phi,\xi,\eta,g)$ be a (2$n$+1) dimensional $\phi$-conharmonically flat Sasakian manifold.\\
Now it is obvious that $\phi^2H(\phi X,\phi Y)\phi Z = 0$ iff
\begin{eqnarray}
  g(H(\phi X,\phi Y)\phi Z,\phi W) &=&0,
\end{eqnarray}
for any $X, Y, Z, W\in\chi(M)$.\\
As the manifold is a (2$n$+1) dimensional $\phi$-conharmonically flat, then using (1.11) and (4.30) we get,
\begin{eqnarray}
  g(R(\phi X,\phi Y)\phi Z,\phi W)&=& \frac{1}{2n-1}[g(\phi Y,\phi Z)S(\phi X,\phi W)\nonumber \\
                                  &-&g(\phi X,\phi Z)S(\phi Y,\phi W)+g(\phi X,\phi W)S(\phi Y,\phi Z)\nonumber\\
                                  &-&g(\phi Y,\phi W)S(\phi X,\phi Z)].
\end{eqnarray}
Let $\{e_1, e_2,.... e_{2n}, \xi\}$ be a local orthonormal basis of $T_p(M)$.\\
Putting $X=W=e_i$ and summing over $i=1,2,....(2n+1)$, we get,
\begin{eqnarray}
 g(R(\phi e_i,\phi Y)\phi Z,\phi e_i)&=& \frac{1}{2n-1}[g(\phi Y,\phi Z)S(\phi e_i,\phi e_i)\nonumber \\
                                  &-&g(\phi e_i,\phi Z)S(\phi Y,\phi e_i)+g(\phi e_i,\phi e_i)S(\phi Y,\phi Z)\nonumber\\
                                  &-&g(\phi Y,\phi e_i)S(\phi e_i,\phi Z)].
\end{eqnarray}
Now in a (2$n$+1) dimensional Sasakian manifold, we have,
\begin{equation}
   g(\phi e_i,\phi e_i)=2n,
\end{equation}
\begin{equation}
   S(\phi e_i,\phi e_i)=r-2n,
\end{equation}
\begin{equation}
 g(\phi e_i,\phi Z)g(\phi Y,\phi e_i)=g(\phi Y,\phi Z),
\end{equation}
\begin{equation}
  g(\phi e_i,\phi Z)S(\phi Y,\phi e_i)=S(\phi Y,\phi Z),
\end{equation}
and
\begin{equation}
 g(R(\phi e_i,\phi Y)\phi Z,\phi e_i)= S(\phi Y,\phi Z)-g(\phi Y,\phi Z).
\end{equation}
Now using (4.33)-(4.37) in (4.32), we get,
\begin{eqnarray}
  S(\phi Y,\phi Z)-g(\phi Y,\phi Z)&=& \frac{1}{2n-1}[(r-2n)g(\phi Y,\phi Z)-S(\phi Y,\phi Z)\nonumber \\
                                  &+&2n S(\phi Y,\phi Z)-S(\phi Y,\phi Z)].
\end{eqnarray}
which implies that,
\begin{equation}
  S(\phi Y,\phi Z)=(r-1)g(\phi Y,\phi Z).
\end{equation}
Now from (3.4), we have,
\begin{equation}
  S(\phi Y,\phi Z)=[2n-1-\lambda+\frac{1}{2}(p + \frac{2}{2n+1})]g(\phi Y,\phi Z).
\end{equation}
Using (3.7) the above equation (4.40) becomes,
\begin{equation}
  S(\phi Y,\phi Z)=[2n+\mu-1]g(\phi Y,\phi Z).
\end{equation}
Again from (3.4), we have,
\begin{equation}
  r=\sum_{i=1}^{2n+1} S(e_i,e_i)=(2n+1)[2n-1-\lambda+\frac{1}{2}(p + \frac{2}{2n+1})]-\mu+1.
\end{equation}
Using (3.7) the above equation (4.42) becomes,
\begin{equation}
  r=2n(2n+\mu).
\end{equation}
Then using (4.41) and (4.43), (4.39) reduces to,
\begin{equation}
  [2n+\mu-1]g(\phi Y,\phi Z)=[2n(2n+\mu)-1]g(\phi Y,\phi Z).
\end{equation}
for any $Y, Z\in \chi(M)$,
and it follows that,
\begin{equation}
  [(2n+\mu-1)-(2n(2n+\mu)-1)]=0,
\end{equation}
which implies that,
\begin{equation}
  \mu=-2n.
\end{equation}
Then using (3.7), we get,
\begin{equation}
  \lambda=\frac{1}{2}(p + \frac{2}{2n+1})+2n.
\end{equation}
This leads to the following:\\
\textbf{Theorem 4.7.} {\em If a (2$n$+1) dimensional Sasakian manifold $ (M,\phi,\xi,\eta,g)$ is $\phi$-conharmonically flat and admits $*$-Conformal $\eta$-Ricci Soliton then $ \mu=-2n$ and $\lambda=\frac{1}{2}(p + \frac{2}{2n+1})+2n$.}\\\\
Now using (4.46), (4.47) and (3.4), we obtain,
\begin{equation}
S(X,Y)=-g(X,Y)+(2n+1)\eta(X)\eta(Y).
\end{equation}
Then we have,\\
\textbf{Corollary 4.8.} {\em If a (2$n$+1) dimensional Sasakian manifold $ (M,\phi,\xi,\eta,g)$ is $\phi$-conharmonically flat and admits $*$-Conformal $\eta$-Ricci Soliton then the manifold becomes $\eta$-Einstein manifold.}\\\\
\textbf{Definition 4.9.} A differentiable manifold $(M^n,g), n>3$, satisfying the condition
\begin{eqnarray}
  \phi^2P(\phi X,\phi Y)\phi Z &=& 0 \nonumber
\end{eqnarray}
is called $\phi$-projectively flat \cite{singh}, where $P$ is the Projective curvature tensor of $M^n$.\\
Let $ (M,\phi,\xi,\eta,g)$ be a (2$n$+1) dimensional $\phi$-projectively flat Sasakian manifold.\\
Now it is obvious that $\phi^2P(\phi X,\phi Y)\phi Z = 0$ iff
\begin{eqnarray}
  g(P(\phi X,\phi Y)\phi Z,\phi W) &=&0,
\end{eqnarray}
for any $X, Y, Z, W\in\chi(M)$.\\
As the manifold is a (2$n$+1) dimensional $\phi$-projectively flat, then using (1.12) and (4.49) we get,
\begin{eqnarray}
  g(R(\phi X,\phi Y)\phi Z,\phi W)&=& \frac{1}{2n}[g(\phi Y,\phi Z)S(\phi X,\phi W)\nonumber\\
                                  &-& g(\phi X,\phi Z)S(\phi Y,\phi W)].
\end{eqnarray}
Let $\{e_1, e_2,.... e_{2n}, \xi\}$ be a local orthonormal basis of $T_p(M)$.\\
Putting $X=W=e_i$ and summing over $i=1,2,....(2n+1)$, we get,
\begin{eqnarray}
 g(R(\phi e_i,\phi Y)\phi Z,\phi e_i)&=& \frac{1}{2n} [g(\phi Y,\phi Z)S(\phi e_i,\phi e_i)\nonumber\\
                                                      &-& g(\phi e_i,\phi Z)S(\phi Y,\phi e_i)].
\end{eqnarray}
Now using (4.33)-(4.37), the above equation (4.51) becomes,
\begin{equation}
  S(\phi Y,\phi Z)-g(\phi Y,\phi Z)=\frac{1}{2n}[(r-2n)g(\phi Y,\phi Z)-S(\phi Y,\phi Z)],
\end{equation}
which implies that,
\begin{equation}
  S(\phi Y,\phi Z)=\frac{r}{2n+1}g(\phi Y,\phi Z).
\end{equation}
Using (4.41) and (4.43) in the above equation (4.53), we get,
\begin{equation}
  \mu=1.
\end{equation}
Then from (3.7), we obtain,
\begin{equation}
  \lambda=\frac{1}{2}(p + \frac{2}{2n+1})-1.
\end{equation}
Then we can state the following theorem:\\
\textbf{Theorem 4.10.} {\em If a (2$n$+1) dimensional Sasakian manifold $ (M,\phi,\xi,\eta,g)$ is $\phi$-projectively flat and admits $*$-Conformal $\eta$-Ricci Soliton then $ \mu=1$ and $\lambda=\frac{1}{2}(p + \frac{2}{2n+1})-1$.}\\\\
Now using (4.54), (4.55) and (3.4), we obtain,
\begin{equation}
S(X,Y)=2ng(X,Y).
\end{equation}
Then we have,\\
\textbf{Corollary 4.11.} {\em If a (2$n$+1) dimensional Sasakian manifold $ (M,\phi,\xi,\eta,g)$ is $\phi$-projectively flat and admits $*$-Conformal $\eta$-Ricci Soliton then the manifold becomes Einstein manifold.}\\\\
\section{\textbf{Example of a 5-dimensional Sasakian manifold:}}
We consider the 5-dimensional manifold $M = \{(x, y, z,u,v) \in \mathbb{R}^5 \}$,
where $(x, y, z, u, v)$ are standard coordinates in $\mathbb{R}^5$. Let ${e_1, e_2, e_3, e_4, e_5}$ be a linearly independent
frame field on M given by,
\begin{equation}
   e_1=2(y \frac{\partial}{\partial z}-\frac{\partial}{\partial x}),\quad e_2 = \frac{\partial}{\partial y}, \quad e_3 =-2\frac{\partial}{\partial z}, \quad e_4=2(v \frac{\partial}{\partial z}-\frac{\partial}{\partial u}), \quad e_5=-2\frac{\partial}{\partial v}. \nonumber
\end{equation}
Let $g$ be the Riemannian metric defined by,
\begin{equation}
  g(e_i,e_j)=0,i\neq j,i,j=1,2,3,4,5,\nonumber
\end{equation}
\begin{equation}
   g(e_1,e_1) = g(e_2,e_2) = g(e_3,e_3) =g(e_4,e_4) = g(e_5,e_5)= 1.\nonumber
\end{equation}
Let $\eta$ be the 1-form defined by $\eta(Z) = g(Z,e_3)$, for any $Z \in \chi(M)$,where $\chi(M)$ is the set of all differentiable vector fields on $M$ and $\phi$ be the (1, 1)-tensor field defined by,
\begin{equation}
  \phi e_1=e_2, \quad \phi e_2=-e_1,\quad  \phi e_3=0, \quad \phi e_4= e_5, \quad \phi e_5=-e_4.\nonumber
\end{equation}
Then, using the linearity of $\phi$ and $g$, we have $\eta(e_3) = 1, \phi ^2 (Z) = -Z + \eta(Z)e_3$ and $g(\phi Z,\phi U) = g(Z,U) - \eta(Z)\eta(U)$, for any $Z,U \in \chi(M)$. Thus for $e_3 = \xi$, $(\phi, \xi, \eta, g)$ defines a Sasakian
structure on $M$.\\
Let $\nabla$ be the Levi-Civita connection with respect to the Riemannian metric $g$. Then we have,\\
  $ [e_1,e_2] =2e_3, [e_4,e_5] =2e_3$, and $ [e_i,e_j] =0$ for others $i,j$.\\
The connection $\nabla$ of the metric $g$ is given by,
\begin{eqnarray}
  2g(\nabla_X Y,Z) &=& Xg(Y,Z)+Yg(Z,X)-Zg(X,Y)\nonumber \\
                   &-& g(X, [Y,Z])-g(Y, [X, Z]) + g(Z, [X, Y ]),\nonumber
\end{eqnarray}
which is known as Koszul’s formula.\\
Using Koszul’s formula, we can easily calculate,
$$\nabla_{e_1} e_1 =0, \nabla_{e_1} e_2 =e_3 , \nabla_{e_1} e_3 =-e_2 , \nabla_{e_1} e_4=0 ,\nabla_{e_1} e_5=0 ,$$
$$\nabla_{e_2} e_1 =-e_3, \nabla_{e_2} e_2 =0, \nabla_{e_2} e_3 =e_1,\nabla_{e_2} e_4=0,\nabla_{e_2} e_5=0,$$
$$\nabla_{e_3} e_1 =-e_2, \nabla_{e_3} e_2 = e_1, \nabla_{e_3} e_3 =0,\nabla_{e_3} e_4=0,\nabla_{e_3} e_5=e_4,$$
$$\nabla_{e_4} e_1=0,\nabla_{e_4} e_2=0,\nabla_{e_4} e_3=-e_5,\nabla_{e_4} e_4=0,\nabla_{e_4} e_5=e_3,$$
$$\nabla_{e_5} e_1=\nabla_{e_5} e_2=\nabla_{e_5} e_3=\nabla_{e_5} e_4=\nabla_{e_5} e_5=0.$$
It can be easily seen that for $e_3 = \xi, (\phi, \xi, \eta, g)$ is a Sasakian structure on $M$. Consequently, $(M, \phi, \xi, \eta, g)$ is a Sasakian manifold.\\
Also, the Riemannian curvature tensor $R$ is given by,
$$R(X, Y )Z = \nabla_X\nabla_Y Z - \nabla_Y \nabla_X Z - \nabla_{[X,Y]} Z.$$
Hence,
$$R(e_1,e_2)e_1 =3e_2, R(e_1,e_3)e_1 = -e_3, R(e_2,e_4)e_1 =-e_5, R(e_2,e_5)e_1=e_4,$$
$$R(e_4,e_5)e_1= 2e_2, R(e_1,e_2)e_2 = -e_1, R(e_1,e_4)e_2 =e_5 ,R(e_2,e_3)e_2=-e_3$$
$$R(e_4,e_5)e_2)=-2e_1,R(e_1,e_3)e_3=e_1,R(e_2,e_3)e_1=-e_3,R(e_3,e_4)e_3=-e_4,$$
$$R((e_4,e_5)e_4=2e_5,R(e_1,e_2)e_5=-2e_4,R(e_1,e_4)e_5=e_2,R(e_2,e_4)e_5=e_1,$$
$$R(e_4,e_5)e_5=-2e_4,R(e_1,e_4)e_5=-e_2.$$
Then, the Ricci tensor $S$ is given by,
$$S(e_1,e_1) =-2, S(e_2,e_2) =3, S(e_3,e_3)= S(e_4,e_4)=4, S(e_5,e_5) = -1.$$
From (3.4) we have, $S(e_3,e_3)=4-\lambda-\mu+\frac{1}{2}(p + \frac{2}{5})$.\\
Therefore we get,
$$4-\lambda-\mu+\frac{1}{2}(p + \frac{2}{5})=4,$$
which implies that,
$$\lambda+\mu=\frac{1}{2}(p + \frac{2}{5}).$$
Hence $\lambda$ and $\mu$ satisfies equation (3.7) and so $g$ defines a $*$-Conformal $\eta$-Ricci Soliton on the 5-dimensional Sasakian manifold $M$.

\end{document}